\documentclass{amsart}
\usepackage{amsmath}
\usepackage{amssymb}
\usepackage{amsthm}
\usepackage{amsfonts}
\usepackage{graphicx}
\newtheorem{Theorem}{Theorem}[section]

\newtheorem{Corollary}[Theorem]{Corollary}

\theoremstyle{definition}

\numberwithin{equation}{section}
\title{On rotational surfaces in 3 dimensional de Sitter space with Weingarten condition}

\author[B. Bekta\c s]{Burcu Bekta\c s Demirci}
\address{Fatih Sultan Vak{\i}f University, Hal\.{I}\c{c} Campus, Faculty of Engineering,
Department of Civil Engineering, 34445, Beyo\u{g}lu, Istanbul}
\email{bbektas@fsm.edu.tr}
\begin{document}

\maketitle
\begin{abstract}
In this article, we study spacelike and timelike rotational surfaces in a 3--dimensional de Sitter space $\mathbb{S}^3_1$ which are the
orbit of a regular curve under the action of the orthogonal
transformation of 4--dimensional Minkowski space $\mathbb{E}^4_1$ leaving 
a spacelike, a timelike or a degenerate plane pointwise fixed. 
We determine the profile curve of 
such Weingarten rotational surfaces parameterized by the principal curvature. 
Then, we classify spacelike and timelike Weingarten rotational surface in $\mathbb{S}^3_1$ 
with the principal curvatures $\kappa$ and $\lambda$ satisfying $\kappa=a\lambda+b$ or $\kappa=a\lambda^m$ for 
special cases of constants $a, b$ and $m$.
\end{abstract}

\section{Introduction}
A surface in a 3--dimensional space forms is called a Weingarten surface, which was first introduced by J. Weingarten in \cite{JW}, if the principal curvatures $\kappa$ and $\lambda$ of the surface satisfy a certain relation $W(\kappa,\lambda)=0$ identically.
Since the mean curvature $H$ and the Gaussian curvature $K$ of the surface are determined by the principal curvatures 
$\kappa$ and $\lambda$, the Weingarten surface also have 
a relation $U(K,H)=0$. 
Thus, minimal surfaces, surfaces with constant mean curvature and surfaces with constant Gaussian curvature are known examples of Weingarten surfaces. 

Many geometers studied Weingarten surfaces as a classical topic of differential geometry and they also have been obtained several results such as 
\cite{HW, RE, KS, L1, L2, L3, Li, Y}. 
Although the complete classification of Weingarten surfaces is still a open problem, the existence results of  
important subclasses of surfaces such as  
tubes along curves and cyclic surfaces \cite{L2,L3}, 
ruled surfaces and helicoidal surfaces \cite{K,DK} and  translation surfaces \cite{DGW, LM} are obtained. 

In this context, 
there are some important classification results for 
rotational surfaces in different ambient spaces. 
For a Weingarten condition $aH+bK=c$ with some constants $a,b$ and $c$, A. Barros et al. \cite{Barros} gave a complete description of all rotational Weingarten surfaces in Euclidean 3--sphere $\mathbb{S}^3$
and R. L\'{o}pez and A. P\'{a}mpano \cite{Lopez} 
make classification of rotational surfaces in Euclidean space satisfying a linear relation between principal curvatures. 

In \cite{KS}, W. Kühnel and M. Stelle studied closed rotational surfaces in Euclidean 3--space satisfying $\kappa=f(\lambda)$ for the principal curvatures $\kappa$ and $\lambda$. 
As a particular case, they obtained analytic closed surfaces of genus zero where $f$ is a quadratic polynomial or $f(\lambda)=c\lambda^{2n+1}$.
Moreover,U. Dursun \cite{UD} gave the results for the profile curve of Weingarten rotational surfaces in the hyperbolic 3--space $\mathbb{H}^3(-1)$ parameterized in terms of principal curvature. 

In a 3--dimensional de Sitter space $\mathbb{S}^3_1$, 
there are three types of rotational surfaces, called spherical, hyperbolic and parabolic rotational surfaces.
In \cite{LL}, H. Liu and G. Liu studied Weingarten rotational surfaces in $\mathbb{S}^3_1$ and 
they determined the coordinates of the profile curve of all spacelike and timelike Weingarten rotational surfaces in $\mathbb{S}^3_1$ parameterized by arc length parameter. 
Unfortunately, this is not true due to the fact that explained in Section 3. 

In this work, we study spacelike and timelike rotational surfaces in $\mathbb{S}^3_1$ with the principal curvatures 
$\kappa$ and $\lambda$ satisfying $\kappa=f(\lambda)$ for a continuous function $f$. First, we gave the parameterization
of the profile curve of such Weingarten surfaces in terms of the principal curvature $\lambda$. Then, we examine special cases of constants $a,b$ and $m$ 
when the Weingarten relation is $f(\lambda)=a\lambda+b$ or $f(\lambda)=a\lambda^m$, respectively. 
For particular choices of constants $a,b$ and $m$ in the Weingarten relations, we obtain classifications of all spacelike and timelike rotational surfaces in $\mathbb{S}^3_1$ 
with constant mean curvature in $\mathbb{S}^3_1$
or constant Gaussian curvature in $\mathbb{S}^3_1$.

\section{Preliminaries}
Let $\mathbb{E}^4_1$ be a 4--dimensional Minkowski space with its metric $\tilde{g}$ is given by 
\begin{equation}
    \tilde{g}({\bf x}, {\bf y})=x_1y_1+x_2y_2+x_3y_3-x_4y_4
\end{equation}
where ${\bf x}=(x_1,x_2,x_3,x_4)$ and 
${\bf y}=(y_1,y_2,y_3,y_4)$.
As a hypersurface of $\mathbb{E}^4_1$, 
the 3--dimensional de Sitter space with constant sectional curvature $1$ denoted by $\mathbb{S}^3_1$ 
is defined by
\begin{equation}
    \mathbb{S}^3_1=\{{\bf x}\in{\mathbb{E}^4_1}\;|\;
    \tilde{g}({\bf x},{\bf x})=1\}
\end{equation}
From \cite{LL} and \cite{Carmo},  
we give the definition of rotational surface in a 3--dimensional de Sitter space $\mathbb{S}^3_1$ as follows.

An orthogonal transformation of $\mathbb{E}^4_1$ is 
a linear map that preserves the metric $\tilde{g}$. 
Let denote an $k$--dimensional subspace of $\mathbb{E}^4_1$ passing through the origin by ${\bf P}^k$ 
and the group of 
orthogonal transformations of $\mathbb{E}^4_1$ with positive determinant that leave ${\bf P}^2$ pointwise fixed
by ${\bf O}({\bf P}^2)$. 
Then, ${\bf P}^k$ is pseudo--Riemannian if the restriction
$\tilde{g}|_{{\bf P}^k}$ is a pseudo--Riemannian metric; 
${\bf P}^k$ is Riemannian if $\tilde{g}|_{{\bf P}^k}$ is a Riemannian metric; 
${\bf P}^k$ is degenerate if $\tilde{g}|_{{\bf P}^k}$ is a degenerate quadratic form.  

We choose ${\bf P}^2$ and ${\bf P}^3$ 
such that ${\bf P}^2\subset {\bf P}^3$ and 
${\bf P}^3\cap\mathbb{S}^3_1\neq\emptyset$. 
Let $C$ be a regular curve in ${\bf P}^3\cap\mathbb{S}^3_1$
that does not intersect ${\bf P}^2$.
The orbit of $C$ under the action of ${\bf O}({\bf P}^2)$ is called rotational surface $M$ in $\mathbb{S}^3_1$ generated by 
$C$ around ${\bf P}^2$
if the induced metric of $M$ from $\mathbb{E}^4_1$ is nondegenerate. The rotational surface $M$ in $\mathbb{S}^3_1$ is said to be spherical
(resp., hyperbolic or parabolic) if 
$\tilde{g}|_{\bf{P}^2}$
is a pseudo--Riemannian metric
(resp., a Riemannian metric or a degenerate  quadratic form).

Let ${\bf x}=(x_1,x_2,x_3,x_4)$ and 
${\bf y}=(y_1,y_2,y_3,y_4)$ be vectors in $\mathbb{E}^4_1$. For the spherical and the hyperbolic rotational surface
in $\mathbb{S}^3_1$, 
we have 
\begin{equation}
 \tilde{g}({\bf x}, {\bf y})=x_1y_1+x_2y_2+x_3y_3-x_4y_4   
\end{equation}
and for the parabolic rotational surface in $\mathbb{S}^3_1$, 
we have 
\begin{equation}
     \tilde{g}({\bf x},{\bf y})=x_1y_1+x_2y_2+x_3y_4+x_4y_3.  
\end{equation}

Let $C$ be a regular curve
which is parameterized by arc length parameter 
$u\in J\subset\mathbb{R}$ 
with $\tilde{g}(C', C')=\varepsilon=\pm 1$
where $\prime$ 
denotes the derivative of the profile curve 
$C$ respect to $u$. 
Then, the parameterization of the rotational surfaces in $\mathbb{S}^3_1$ are given as follows.

The spherical rotational surface $M$ in $\mathbb{S}^3_1$ 
is defined by 
\begin{equation}
\label{sphr}
    r(u,v)=(y(u)\sin{v},y(u)\cos{v},z(u),w(u))\;\;
    u\in J,\;v\in[0,2\pi].
\end{equation}
The induced metric $g$ on $M$ is given by 
\begin{equation}
    g=\varepsilon du^2+y^2(u)dv^2.
\end{equation}
It can be easily seen that the surface $M$ is spacelike 
for $\varepsilon=1$ and the surface $M$ is timelike
for $\varepsilon=-1$. 
The coordinate functions $y(u), z(u)$ and $w(u)$ of profile curve $C$ of $M$ satisfy
\begin{equation}
    y^2(u)+z^2(u)-w^2(u)=1\;\;\mbox{and}\;\;
    y'^2(u)+z'^2(u)-w'^2(u)=\varepsilon
\end{equation}
and we assume $y(u)>0$ on $J\subset\mathbb{R}$. 
Since $M$ is not contained in a hyperplane of $\mathbb{E}^4_1$, $z'(u)$ and $w'(u)$ can not be zero. 

The hyperbolic rotation surface $M$ of the first kind in $\mathbb{S}^3_1$ is defined by
\begin{equation}
\label{hyp1}
    r(u,v)=(x(u),y(u),w(u)\sinh{v},w(u)\cosh{v}),\;\;
    u\in J,\; v\in\mathbb{R}.
\end{equation}
The induced metric 
$g$ on $M$ is given by 
\begin{equation}
    g=\varepsilon du^2+w^2(u)dv^2.
\end{equation}
Then, the surface $M$ is spacelike for $\varepsilon=1$ 
and the surface $M$ is timelike for $\varepsilon=-1$. 
The coordinate functions $x(u), y(u)$ and $w(u)$ of the profile curve $C$ of $M$ satisfy
\begin{equation}
   x^2(u)+y^2(u)-w^2(u)=1\;\;\mbox{and}\;\;
    x'^2(u)+y'^2(u)-w'^2(u)=\varepsilon
\end{equation}
and we assume $w(u)>0$ on $J\subset\mathbb{R}$. 
Since $M$ is not contained in a hyperplane of $\mathbb{E}^4_1$, $x'(u)$ and $y'(u)$ can not be zero.

The hyperbolic rotation surface $M$ of the second kind in $\mathbb{S}^3_1$ is defined by
\begin{equation}
\label{hyp2}
    r(u,v)=(x(u),y(u),w(u)\cosh{v},w(u)\sinh{v}),\;
    u\in J,\; v\in\mathbb{R}.
\end{equation}
The induced metric $g$ on $M$ is given by 
\begin{equation}
    g=du^2-w^2(u)dv^2.
\end{equation}
Thus, the surface $M$ is timelike. 
Also, the coordinate functions $x(u), y(u)$ and $w(u)$ of the profile curve $C$ of $M$ satisfy the following equations
\begin{equation}
   x^2(u)+y^2(u)+w^2(u)=1\;\;\mbox{and}\;\;
    x'^2(u)+y'^2(u)+w'^2(u)=1
\end{equation}
and we assume that $w(u)>0$ on $J\subset\mathbb{R}$. 
Since $M$ is not contained in a hyperplane of $\mathbb{E}^4_1$, $x'(u)$ and $y'(u)$ can not be zero.

The parabolic rotation surface $M$ in $\mathbb{S}^3_1$ 
is defined by
\begin{equation}
\label{prb}
    r(u,v)=\left(x(u),vz(u),z(u),-\frac{1}{2}v^2z(u)+w(u)
    \right),\;u\in J,\; v\in\mathbb{R}.
\end{equation}
The induced metric 
$g$ on $M$ is given by 
\begin{equation}
    g=\varepsilon du^2+z^2(u)dv^2.
\end{equation}
When $\varepsilon=1$, the surface $M$ is spacelike; when 
$\varepsilon=-1$, the surface $M$ is timelike. 
The coordinate functions $x(u), z(u)$ and $w(u)$ of the profile curve $C$ of $M$ satisfy
\begin{equation}
   2z(u)w(u)+x^2(u)=1\;\;\mbox{and}\;\;
    2z'(u)w'(u)+x'^2(u)=\varepsilon
\end{equation}
and we assume $z(u)>0$ on $J\subset\mathbb{R}$.
Since $M$ is not contained in a hyperplane of $\mathbb{E}^4_1$, $x'(u)$ and $z'(u)$ can not be zero. 

From now on, we denote a rotational surface $M$ in $\mathbb{S}^3_1$ by a notation 
${\bf M}_\delta$ for $\delta=1,-1,0$. 
According to the values of $\delta$,
${\bf M}_1$(resp., ${\bf M}_{-1}$ or ${\bf M}_0$) is a spherical rotational surface in $\mathbb{S}^3_1$ defined by \eqref{sphr}
(resp., a hyperbolic rotational surface of first kind in $\mathbb{S}^3_1$ defined by \eqref{hyp1} or 
a parabolic rotational surface in $\mathbb{S}^3_1$ defined by \eqref{prb}).

%%%%%%%%%%%%%%%%%%%%%%%%%%%%%%%%%%%%%%%%%%%%%%%%%%%
\section{Weingarten rotational surfaces in $\mathbb{S}^3_1$}
In \cite{LL}, H. Liu and G. Liu examined rotational surfaces in $\mathbb{S}^3_1$ and they gave classification theorem for all spacelike and timelike Weingarten rotational surface in $\mathbb{S}^3_1$. 
In their results, the coordinate functions of the profile curve 
of Weingarten rotational surfaces were determined in terms of arc length parameter which is not true. 
On the other hand, the solutions of the obtained differential equations corresponding to Weingarten condition $\kappa=f(\lambda)$ depend on 
principal curvature $\lambda$ not arc length parameter $u$. 

Thus, we can give the following theorem by using computations in \cite{LL}.

\begin{Theorem}
\label{Weingarten1}
Let ${\bf M}_\delta$ be a Weingarten rotational surface in a 3--dimensional de Sitter space $\mathbb{S}^3_1$ 
with the principal curvatures $\kappa$ and $\lambda$ satisfying $\kappa=f(\lambda)$ 
for some continuous function $f$.
Assume ${\bf M}_\delta$ has no umbilical points. 
Then, the parameterization of the surface ${\bf M}_\delta$ with respect to $\lambda$ 
is given by one of the followings:
\begin{itemize}
    \item [(i.)] the spherical rotational surface 
    ${\bf M}_1$
    \begin{equation}
\label{spherical1}
    {\bf{r}_1}(\lambda,v)=\left(y(\lambda)\sin{v}, y(\lambda)\cos{v},z(\lambda),w(\lambda)\right),
    \;\;\lambda\in J\subset\mathbb{R},\;v\in[0,2\pi]
    \end{equation}
    \begin{equation}
    \label{spherical3}
      z(\lambda)=\sqrt{1-y^2(\lambda)}\cosh{\phi(\lambda)}\;\;\mbox{and}\;\; 
    w(\lambda)=\sqrt{1-y^2(\lambda)}\sinh{\phi(\lambda)}\;\;
    \mbox{for}\;\; 0<y(\lambda)<1  
    \end{equation}
    or 
    \begin{equation}
    \label{spherical4}
    z(\lambda)=\sqrt{y^2(\lambda)-1}\sinh{\phi(\lambda)}\;\;\mbox{and}\;\; 
    w(\lambda)=\sqrt{y^2(\lambda)-1}\cosh{\phi(\lambda)}
    \;
    \mbox{for}\;\; y(\lambda)>1
    \end{equation}
    where 
    \begin{equation}
\label{spherical2}
y(\lambda)=ce^{\int{\frac{d\lambda}{\varepsilon f(\lambda)-\lambda}}}
\end{equation}
and 
\begin{equation}
\label{spherical5}
    \phi(\lambda)=\phi_0\pm\int{\frac{\lambda y^2(\lambda)d\lambda}{(y^2(\lambda)-1)(\varepsilon f(\lambda)-\lambda)
    \sqrt{(\lambda^2-\varepsilon)y^2(\lambda)+\varepsilon}}}
\end{equation}
with constants $c>0$ and $\phi_0$, $\varepsilon=\pm 1$ and $(\lambda^2-\varepsilon)y^2(\lambda)+\varepsilon>0$ on an open interval $\tilde{J}\subset J$.
When $\varepsilon=1$, 
the surface ${\bf M}_1$ is spacelike; 
when 
$\varepsilon=-1$, the surface ${\bf M}_1$ is timelike. 

\item [(ii.)] the hyperbolic rotational surface of first kind ${\bf M}_{-1}$
\begin{equation}
\label{hyper1kind1}
    {\bf{r}_{-1}}(\lambda,v)=\left(x(\lambda), y(\lambda), w(\lambda)\sinh{v}, w(\lambda)\cosh{v}\right),\;\;
    \lambda\in J\subset\mathbb{R},\;v\in\mathbb{R}
\end{equation}
\begin{equation}
     \label{hyper1kind3}
     x(\lambda)=\sqrt{1+w^2(\lambda)}\cos{\phi(\lambda)}\;\;\mbox{and}\;\; 
     y(\lambda)=\sqrt{1+w^2(\lambda)}\sin{\phi(\lambda)}\;\;
     \mbox{for}\;\;w(\lambda)>0
    \end{equation}
where 
\begin{equation}
\label{hyper1kind2}
w(\lambda)=ce^{\int{\frac{d\lambda}{\varepsilon f(\lambda)-\lambda}}}
\end{equation}
\end{itemize}
and 
\begin{equation}
 \label{hyper1kind4}
 \phi(\lambda)=\phi_0\pm\int{\frac{\lambda w^2(\lambda) d\lambda}{(1+w^2(\lambda))(\varepsilon f(\lambda)-\lambda)\sqrt{(\lambda^2-\varepsilon)w^2(\lambda)-\varepsilon}}}
\end{equation}
with constants $c>0$ and $\phi_0$, $\varepsilon=\pm 1$ and 
$(\lambda^2-\varepsilon)w^2(\lambda)-\varepsilon>0$ on an open interval 
$\tilde{J}\subset J$.
When $\varepsilon=1$, the surface ${\bf M}_{-1}$ 
is spacelike; 
when $\varepsilon=-1$, the surface ${\bf M}_{-1}$ 
is timelike. 

\item {(iii.)} the hyperbolic rotational surface of second kind $M$
\begin{equation}
\label{hyper2kind1}
    {\bf{r}}(\lambda,v)=\left(x(\lambda), y(\lambda), w(\lambda)\cosh{v}, w(\lambda)\sinh{v}\right),\;\;
    \lambda\in J\subset\mathbb{R},\;v\in\mathbb{R}
\end{equation}
\begin{equation}
     \label{hyper2kind3}
     x(\lambda)=\sqrt{1-w^2(\lambda)}\cos{\phi(\lambda)}\;\;\mbox{and}\;\; 
     y(\lambda)=\sqrt{1-w^2(\lambda)}\sin{\phi(\lambda)}\;\;
     \mbox{for}\;\;0<w(\lambda)<1
    \end{equation}
where 
\begin{equation}
\label{hyper2kind2}
w(\lambda)=ce^{-\int{\frac{d\lambda}{f(\lambda)+\lambda}}}
\end{equation}
and 
\begin{equation}
 \label{hyper2kind4}
 \phi(\lambda)=\phi_0\pm\int{\frac{\lambda w^2(\lambda) d\lambda}{(1-w^2(\lambda))(f(\lambda)+\lambda)
 \sqrt{1-(\lambda^2+1)w^2(\lambda)}}}
\end{equation}
with constants $c>0$, $\phi_0$ and $1-(\lambda^2+1)w^2(\lambda)>0$ on 
an open interval $\tilde{J}\subset J$. 
Then, the surface $M$ is timelike.

\item {(iv.)} the parabolic rotational surface ${\bf M}_0$
\begin{equation}
\label{parabolic1}
    {\bf{r}_0}(\lambda,v)=\left(x(\lambda),vz(\lambda),z(\lambda),-\frac{1}{2}v^2z(\lambda)+w(\lambda)\right),\;\;\lambda\in J\subset\mathbb{R},\;v\in\mathbb{R}
\end{equation}
\begin{equation}
\label{parabolic3}
x(\lambda)=\phi(\lambda)z(\lambda)\;\;\;\mbox{and}\;\;\;
w(\lambda)=\frac{1-x^2(\lambda)}{2z(\lambda)}\;\;
\mbox{for}\;\;z(\lambda)>0
\end{equation}
where 
\begin{equation}
\label{parabolic2}
z(\lambda)=ce^{\int{\frac{d\lambda}{\varepsilon f(\lambda)-\lambda}}}
\end{equation}
and
\begin{equation}
\label{parabolic4}
\phi(\lambda)=\phi_0\pm\int\frac{\lambda d\lambda}{z(\lambda)
(\varepsilon f(\lambda)-\lambda)\sqrt{\lambda^2-\varepsilon}}
\end{equation}
with constants $c>0$ and $\phi_0$, $\varepsilon=\pm 1$
and $\lambda^2-\varepsilon>0$ on an open interval $\tilde{J}\subset J$. 
When $\varepsilon=1$, the surface ${\bf M}_0$ is spacelike; when $\varepsilon=-1$, the surface ${\bf M}_0$ is timelike.
\end{Theorem}

\subsection{Weingarten rotational surfaces in $\mathbb{S}^3_1$ with $f(\lambda)=a\lambda+b$}
Suppose that the principal curvatures of $\kappa$ and $\lambda$ of the rotational surfaces in $\mathbb{S}^3_1$ satisfy $\kappa=f(\lambda)=a\lambda+b$ where $a$ and $b$ are constants. 

For the spherical rotational surface ${\bf M}_1$ in 
$\mathbb{S}^3_1$ defined by \eqref{spherical1}, 
the equation \eqref{spherical2} gives
\begin{equation}
    \label{spherical6}
    y(\lambda)=c((\varepsilon a-1)\lambda+\varepsilon b)^{1/{\varepsilon a-1}}
\end{equation}
with constants $c>0$, $a\neq\varepsilon$ and 
$(\varepsilon a-1)\lambda+\varepsilon b>0$.

Similarly, the equations \eqref{hyper1kind2} and \eqref{parabolic2}
also give the components $w(\lambda)$ and $z(\lambda)$ of the profile curve of 
the hyperbolic rotational surface of first kind ${\bf M}_{-1}$ 
and parabolic rotational surface 
${\bf M}_0$ in $\mathbb{S}^3_1$ defined by \eqref{hyper1kind1} and \eqref{parabolic1}, respectively, as
\begin{equation}
    \label{spherical6a}
    w(\lambda)=z(\lambda)=c((\varepsilon a-1)\lambda+\varepsilon b)^{1/{\varepsilon a-1}}
\end{equation}
with constants $c>0$, $a\neq\varepsilon$ and 
$(\varepsilon a-1)\lambda+\varepsilon b>0$. 

Let define the functions $W(\lambda)$ and $\Phi_{\varepsilon}(\lambda,\delta)$ as follows 
\begin{align}
\label{eq}
W(\lambda)&=(\varepsilon a-1)\lambda+\varepsilon b,\\
  \label{eq1}
   \Phi_{\varepsilon}(\lambda,\delta)&=\phi_0\pm c^2
    \int{\frac{\lambda W(\lambda)
    ^{\frac{3-\varepsilon a}{\varepsilon a-1}}d\lambda}
    {(c^2 W(\lambda)^{2/ \varepsilon a-1}-\delta)
    \sqrt{c^2(\lambda^2-\varepsilon) W(\lambda)^{2/{\varepsilon a-1}}+\delta\varepsilon}}}. 
\end{align}
Here, $\lambda$ is in some open interval $\tilde{J}\subset J$ on which  
$(\varepsilon a-1)\lambda+\varepsilon b>0$,
 $c^2W(\lambda)^{2/\varepsilon a-1}-\delta\neq 0$
 and
$c^2(\lambda^2-\varepsilon)W(\lambda)^{2/{\varepsilon a-1}}+\delta\varepsilon>0$ are satisfied. 
Then, we can express the components of the profile curve 
of the surfaces ${\bf M}_{\delta}$ in \eqref{spherical6} and 
\eqref{spherical6a} as
$y(\lambda)=w(\lambda)=z(\lambda)=cW(\lambda)
^{1/{\varepsilon a-1}}$.

Also, it can be seen that  for the surfaces ${\bf M}_\delta$,
the function $\phi(\lambda)$ given by \eqref{spherical5}, \eqref{hyper1kind4} and \eqref{parabolic4} can be obtained from the equation \eqref{eq1} substituting the values of $\delta=1,-1,0$, respectively. 
Thus, the remaining coordinate functions of the profile curve for the surfaces ${\bf M}_{\delta}$
in \eqref{spherical3} or \eqref{spherical4}, \eqref{hyper1kind3} and \eqref{parabolic3} are completely determined.  

Now, we will consider case for constants $a=-\varepsilon$ 
and $b=2\varepsilon H$ with constant $H$. 
Using Theorem \ref{Weingarten1}, we give following classification results 
of spacelike and timelike rotational surfaces with constant mean curvature in a 3--dimensional de Sitter space $\mathbb{S}^3_1$ whose profile curve is 
parameterized with respect to the principal curvature $\lambda$. 

\begin{Corollary}
\label{spacelikeWeingarten}
Let ${\bf M}_\delta$ be a spacelike Weingarten rotational surface in a 3--dimensional de Sitter space $\mathbb{S}^3_1$
defined by \eqref{spherical1}, \eqref{hyper1kind1} and \eqref{parabolic1} 
with the principal curvatures $\kappa$ and $\lambda$ satisfying $\kappa=-\lambda+2H$ and 
 \begin{equation}
        \Phi_1(\lambda,\delta)=
        \phi_0\pm\frac{c^2}{\sqrt{2}}\int{\frac{\lambda d\lambda}{\sqrt{H-\lambda}(c^2-2\delta(H-\lambda))\sqrt{(\lambda^2-1)c^2+2\delta(H-\lambda)}}}
    \end{equation}
with constants $\phi_0$, $H$ and $c>0$. 
Then, the coordinate functions of the profile curve of ${\bf M}_\delta$ with constant mean curvature $H$ parameterized with respect to $\lambda$
are given by one of the following cases:
\begin{itemize}
    \item [1.] for $\delta=1$, 
    \subitem (a.) $y(\lambda)=\frac{c}{\sqrt{2(H-\lambda)}},\; z(\lambda)=\sqrt{1-y^2(\lambda)}\cosh{\Phi_1(\lambda,1)},\; w(\lambda)=\sqrt{1-y^2(\lambda)}
    \sinh{\Phi_1(\lambda,1)}$
    with $\lambda<H-\frac{c^2}{2}$,
   
    \subitem (b.) $y(\lambda)=\frac{c}{\sqrt{2(H-\lambda)}},\; z(\lambda)=\sqrt{y^2(\lambda)-1}\sinh{\Phi_1(\lambda,1)},\; w(\lambda)=\sqrt{y^2(\lambda)-1}
    \cosh{\Phi_1(\lambda,1)}$
    where
    \subsubitem i)
    $H-\frac{c^2}{2}<\lambda< H$ when $H\geq\frac{1}{2}\left(c^2+\frac{1}{c^2}\right)$ or 
    $1\leq H<\frac{1}{2}\left(c^2+\frac{1}{c^2}\right)$ for $c<1$
    or $H\leq  -1$,
    
    \subsubitem ii)
    $H-\frac{c^2}{2}<\lambda
    <\frac{1-\sqrt{1-2Hc^2+c^4}}{c^2}$ when $H<\frac{c^2}{2}$ or $\frac{c^2}{2}<H
    <\frac{1}{2}\left(c^2+\frac{1}{c^2}\right)$
    
    \subsubitem iii) $\frac{1+\sqrt{1-2Hc^2+c^4}}{c^2}<\lambda<H$ when $1<H<\frac{1}{2}\left(c^2+\frac{1}{c^2}\right)$ for $c>1$,
    
    \item [2.] for $\delta=-1$, 
    $$w(\lambda)=\frac{c}{\sqrt{2(H-\lambda)}},\; x(\lambda)=\sqrt{1+w^2(\lambda)}\cos{\Phi_1(\lambda,-1)},
    \; 
    y(\lambda)=\sqrt{1+w^2(\lambda)}\sin{\Phi_1(\lambda,-1)}$$
    where 
    \subitem i) $\lambda<H$ when $H\leq-\frac{1}{c^2}$ for $c>1$,
    \subitem ii) $\displaystyle{\lambda<\frac{-1-\sqrt{1+2Hc^2+c^4}}{c^2}}$ when
    $-1\leq H\leq-\frac{1}{c^2}$ for $c\geq 1$,
    \subitem iii) 
    $\displaystyle{\frac{-1+\sqrt{1+2Hc^2+c^4}}{c^2}}
    <\lambda<H$ when
    $-\frac{1}{2}\left(c^2+\frac{1}{c^2}\right)<H\leq -1$ for $c<1$ or
    $H\geq 1$ for $c>1$.

    \item [3.] for $\delta=0$,
    $$z(\lambda)=\frac{c}{\sqrt{2(H-\lambda)}},\; x(\lambda)=\Phi_1(\lambda,0)z(\lambda),\; w(\lambda)=\frac{1-x^2(\lambda)}{2z(\lambda)}$$
    where 
    \subitem i) $\lambda<H$ when $H<-1$,
    \subitem ii) $\lambda<-1$ when $-1\leq H \leq 1$,
    \subitem iii) $\lambda<-1$ or $1<\lambda<H$ when $H>1$. 
\end{itemize}
\end{Corollary}

\begin{Corollary}
\label{timelikeWeingarten}
Let ${\bf M}_\delta$ be a timelike Weingarten rotational surface in a 3--dimensional de Sitter space $\mathbb{S}^3_1$
defined by \eqref{spherical1}, \eqref{hyper1kind1} and \eqref{parabolic1} 
with the principal curvatures $\kappa$ and $\lambda$
satisfying $\kappa=\lambda-2H$ and 
\begin{equation}
        \Phi_{-1}(\lambda,\delta)=\phi_0\pm\frac{c^2}{\sqrt{2}}\int{\frac{\lambda d\lambda}{\sqrt{H-\lambda}(c^2-2\delta(H-\lambda))\sqrt{(\lambda^2+1)c^2-2\delta(H-\lambda)}}}
    \end{equation}
with constants $\phi_0$, $H$ and $c>0$. 
Then, the coordinate functions of the profile curve of ${\bf M}_\delta$ with the constant mean curvature $H$ in $\mathbb{S}^3_1$ with respect to $\lambda$
are given by one of the following cases:
\begin{itemize}
    \item [1.] for $\delta=1$,
    \subitem (a.) $y(\lambda)=\frac{c}{\sqrt{2(H-\lambda)}},\; z(\lambda)=\sqrt{1-y^2(\lambda)}
    \cosh{\Phi_{-1}(\lambda,1)},\; w(\lambda)=\sqrt{1-y^2(\lambda)}
    \sinh{\Phi_{-1}(\lambda,1)}$
    where
    \subsubitem i) $\lambda<\frac{-1-\sqrt{1+2Hc^2-c^4}}{c^2}$ or 
    $\frac{-1+\sqrt{1-2Hc^2-c^4}}{c^2}<
    \lambda<H-\frac{c^2}{2}$ when $H>\frac{1}{2}\left(c^2-\frac{1}{c^2}\right)$
    \subsubitem ii) 
    $\lambda<H-\frac{c^2}{2}$ when $H\leq\frac{1}{2}\left(c^2-\frac{1}{c^2}\right)$

    \subitem (b.) $y(\lambda)=\frac{c}{\sqrt{2(H-\lambda)}},\; z(\lambda)=\sqrt{y^2(\lambda)-1}
    \sinh{\Phi_{-1}(\lambda,1)},\; w(\lambda)=\sqrt{y^2(\lambda)-1}
    \cosh{\Phi_{-1}(\lambda,1)}$
    where $H-\frac{c^2}{2}<\lambda<H$,
    
    \item [2.] for $\delta=-1$, 
    $$w(\lambda)=\frac{c}{\sqrt{2(H-\lambda)}},\; x(\lambda)=\sqrt{1+w^2(\lambda)}
    \cos{\Phi_{-1}(\lambda,-1)},\; y(\lambda)=\sqrt{1+w^2(\lambda)}
    \sin{\Phi_{-1}(\lambda,-1)}$$
    where $\lambda<H$,
    
    \item [3.] for $\delta=0$,
    $$z(\lambda)=\frac{c}{\sqrt{2(H-\lambda)}},\; x(\lambda)=\Phi_{-1}(\lambda,0)z(\lambda),\; w(\lambda)=\frac{1-x^2(\lambda)}{2z(\lambda)}$$
    where $\lambda<H$.  
\end{itemize}
\end{Corollary}

As doing similar calculation for the hyperbolic rotational surface of second kind $M$ in $\mathbb{S}^3_1$ defined by \eqref{hyper2kind1}
the equation \eqref{hyper2kind1} gives the following
\begin{equation}
    \label{hyper2kind5}
    w(\lambda)=c((a+1)\lambda+b)^{-1/{a+1}}
\end{equation}
with constants $c>0$, $a\neq -1$ and
$(a+1)\lambda+b>0$. 
Let say $\overline{W}(\lambda)=(a+1)\lambda+b$. 
Then, we get $w(\lambda)=c\overline{W}(\lambda)^{-1/{a+1}}$. 
Thus, the equation \eqref{hyper2kind4} gives following equation:
\begin{equation}
    \label{hyper2kind6}
    \phi(\lambda)=\phi_0\pm c^2
    \int{\frac{\lambda \overline{W}(\lambda)^{-\frac{3+a}{a+1}}d\lambda}{
    (1-c^2 \overline{W}(\lambda)^{-2/{a+1}})
    \sqrt{1-(\lambda^2+1)c^2\overline{W}(\lambda)^{-2/{a+1}}}}}.
\end{equation}
where $\lambda$ is in some open interval $\tilde{J}\subset J$ on which
$(a+1)\lambda+b>0$, $1-c^2\bar{W}(\lambda)^{-2/a+1}\neq 0$ 
and 
$1-c^2(\lambda^2+1)\bar{W}(\lambda)^{-2/{a+1}}>0$ are satisfied.
Thus, the remaining coordinate functions of the profile curve of $M$
are found from \eqref{hyper2kind3}.

Now, we will take constants $a$ and $b$ as $a=1$ and $b=2H$ with constant $H$. 
From Theorem \ref{Weingarten1}, we give
the parametrization of hyperbolic rotational surfaces of second kind with constant mean curvature in a 3--dimensional de Sitter space $\mathbb{S}^3_1$ 
in terms of the principal curvature $\lambda$ as below.
\begin{Corollary}
\label{hyper2kindcor1}
Let $M$ be a Weingarten hyperbolic rotational of second kind in a 3--dimensional de Sitter space $\mathbb{S}^3_1$ defined by \eqref{hyper2kind1} with the principal curvatures 
$\kappa$ and
$\lambda$ satisfying $\kappa=\lambda+2H$
with constant $H$. Then, the coordinate functions of the profile curve
of $M$ with constant mean curvature $H$ in $\mathbb{S}^3_1$ parameterized with respect to $\lambda$
are given as follows:
    \begin{equation}
        \label{hyper2kind7}
        w(\lambda)=\frac{c}{\sqrt{2(H+\lambda)}},\;\; 
        x(\lambda)=\sqrt{1-w^2(\lambda)}\cos{\phi(\lambda)}\;\;
        y(\lambda)=\sqrt{1-w^2(\lambda)}\sin{\phi(\lambda)}
    \end{equation}
where 
\begin{equation}
    \label{hyper2kind8}
    \phi(\lambda)=\phi_0\pm\frac{c^2}{\sqrt{2}}
    \int{
    \frac{\lambda d\lambda}
    {\sqrt{H+\lambda}(2(H+\lambda)-c^2)\sqrt{2(H+\lambda)-(\lambda^2+1)c^2}}}
\end{equation}
with 
$\frac{1-\sqrt{1+2Hc^2-c^4}}{c^2}<\lambda<\frac{1+\sqrt{1+2Hc^2-C^4}}{c^2}$
provided that 
$H>\frac{1}{2}\left(c^2-\frac{1}{c^2}\right)$ for 
constants $\phi_0$ and $c>0$.
\end{Corollary}

\subsection{Weingarten rotational surface in $\mathbb{S}^3_1$ with $f(\lambda)=a\lambda^m$}
Assume that the principal curvatures $\kappa$ and $\lambda$
of the rotational surfaces in $\mathbb{S}^3_1$ satisfy 
$\kappa=f(\lambda)=a\lambda^m$ for constants $a$ and $m\neq 1$. 

For the spherical rotational surface ${\bf M}_1$ in $\mathbb{S}^3_1$ defined by \eqref{spherical1}, 
the equation \eqref{spherical2} implies 
\begin{equation}
  \frac{dy}{y}=
  \frac{d\lambda}{\varepsilon f(\lambda)-\lambda}
\end{equation}
and if we substitute $f(\lambda)=a\lambda^m$ into this equation,
then we get the following differential equation
\begin{equation}
    \frac{d\lambda}{dy}+\frac{\lambda}{y}=\frac{\varepsilon a}{y}\lambda^m.
\end{equation}
which is called  Bernoulli equation for $m\neq 1$. 
Then, its solution is 
\begin{equation}
    \label{spherical16}
    y(\lambda)=c(\lambda^{1-m}-\varepsilon a)^{1/{m-1}}
\end{equation}
with a positive constant $c$. 

From similar calculation, we get the components $w(\lambda)$ and $z(\lambda)$
in \eqref{hyper1kind1} and \eqref{parabolic1}
\begin{equation}
    \label{spherical161}
    w(\lambda)=z(\lambda)=c(\lambda^{1-m}-\varepsilon a)^{1/{m-1}}.
\end{equation}
These are the coordinate functions of the profile curve of 
the hyperbolic rotational surface of first kind ${\bf M}_{-1}$ and 
the parabolic rotational surface ${\bf M}_0$ defined by \eqref{hyper1kind1} 
and \eqref{parabolic1}, respectively.

Let define the functions $S(\lambda)$ and $\Psi_{\varepsilon}(\lambda,\delta)$
as 
\begin{align}
\label{eqA}
S(\lambda)&=\lambda^{1-m}-\varepsilon a\\
    \label{eq2}
    \Psi_{\varepsilon}(\lambda,\delta)&=\psi_0\pm c^2\int{\frac{\lambda^{1-m}S(\lambda)^{\frac{3-m}{m-1}}d\lambda}
    {(c^2 S(\lambda)^{2/{m-1}}-\delta)\sqrt{(\lambda^2-\varepsilon)c^2
    S(\lambda)^{2/{m-1}}+\delta\varepsilon}}}.
\end{align}
Then, we have the coordinate functions of the surface ${\bf M}_\delta$ as 
$y(\lambda)=w(\lambda)=z(\lambda)=cS(\lambda)^{1/m-1}$.

Moreover, the function $\Psi_{\varepsilon}(\lambda,\delta)$
defined in \eqref{eq2} gives the function $\phi(\lambda)$ for the Weingarten rotational surfaces ${\bf M}_{\delta}$ according to the values of $\delta=1,-1,0$, respectively. 
Thus, the remaining coordinate functions of the profile curve of 
${\bf M}_{\delta}$
are given by \eqref{spherical3} or \eqref{spherical4}, \eqref{hyper1kind3} and \eqref{parabolic3}.

Now, we will consider case for $m=-1$. 
From Theorem \ref{Weingarten1}, we get following classification results of 
the spacelike and timelike rotational surfaces with constant Gaussian curvature in $\mathbb{S}^3_1$ whose profile curve is parameterized by the principal curvature $\lambda$. 
\begin{Corollary}
\label{spacelikeWeingarten1}
Let ${\bf M}_\delta$ be a spacelike Weingarten rotational surfaces in a 3--dimensional de Sitter space $\mathbb{S}^3_1$ defined by \eqref{spherical1}, \eqref{hyper1kind1} and 
\eqref{parabolic1} with the principal curvatures $\kappa$ and $\lambda$
satisfying $\kappa=\frac{a}{\lambda}$ and 
\begin{equation}
    \Psi_{1}(\lambda,\delta)=\psi_0\pm 
    c^2\int{\frac{\lambda^2 d\lambda}{\sqrt{\lambda^2-a}(c^2-\delta(\lambda^2-a))
    \sqrt{(\lambda^2-1)c^2+\delta(\lambda^2-a)}}}
\end{equation}
for constants $\psi_0$, $a$ and $c>0$.
Then, the coordinate functions of the profile curve of ${\bf M}_\delta$ with 
constant Gaussian curvature in $\mathbb{S}^3_1$ parameterized with respect to 
$\lambda$ are given by one of the following:
\begin{itemize}
    \item [1.] for $\delta=1$,
    \subitem (a.)
    $y(\lambda)=\frac{c}{\sqrt{\lambda^2-a}},\; 
    z(\lambda)=\sqrt{1-y^2(\lambda)}\cosh{\Psi_1(\lambda,1)},\;
    w(\lambda)=\sqrt{1-y^2(\lambda)}\sinh{\Psi_1(\lambda,1)}$
    where
    \subsubitem i) $\lambda<-\sqrt{c^2+a}$ or $\lambda>\sqrt{c^2+a}$ for $a>-c^2$,
    \subsubitem ii) $\lambda\in\mathbb{R}$ for 
    $a\leq -c^2$,
    
    \subitem (b.) 
     $y(\lambda)=\frac{c}{\sqrt{\lambda^2-a}},\; 
        z(\lambda)=\sqrt{y^2(\lambda)-1}
        \sinh{\Psi_1(\lambda,1)},\;
        w(\lambda)=\sqrt{y^2(\lambda)-1}
        \cosh{\Psi_1(\lambda,1)}$
        where
        \subsubitem i) $-\sqrt{c^2+a}<\lambda<-\sqrt{\frac{c^2+a}{c^2+1}}$ or $\sqrt{\frac{c^2+a}{c^2+1}}<\lambda<\sqrt{c^2+a}$ for $-c^2<a<1$,
        \subsubitem i) $-\sqrt{c^2+a}<\lambda<-\sqrt{a}$ or $\sqrt{a}<\lambda<\sqrt{c^2+a}$ for $a>1$,
        
    \item [2.] for $\delta=-1$,
        $$w(\lambda)=\frac{c}{\sqrt{\lambda^2-a}},\;
        x(\lambda)=\sqrt{1+w^2(\lambda)}
        \cos{\Psi_1(\lambda,-1)},
        \;\;
        y(\lambda)=\sqrt{1+w^2(\lambda)}
        \sin{\Psi_1(\lambda,-1)}$$
       where 
       \subitem i) $\lambda<-\sqrt{a}$ or $\lambda>\sqrt{a}$ for $a>1$ and $c>1$,
       
       \subitem ii) $\lambda<-\sqrt{\frac{c^2-a}{c^2-1}}$ or $\lambda>\sqrt{\frac{c^2-a}{c^2-1}}$ for $a<1<c^2$, 
       
     \item [3.] for $\delta=0$, 
      $$z(\lambda)=\frac{c}{\sqrt{\lambda^2-a}},\;
        x(\lambda)=\Psi_{1}(\lambda,0)z(\lambda),
        \;\;
        w(\lambda)=\frac{1-x^2(\lambda)}{2z(\lambda)}$$
       with 
       \subitem i) $\lambda<-\sqrt{a}$ or $\lambda>\sqrt{a}$ for $a>1$,
       
       \subitem ii) $\lambda<-1$ or $\lambda>1$ for $a\leq 1$.
\end{itemize}
\end{Corollary}

\begin{Corollary}
\label{timelikeWeingarten1}
Let ${\bf M}_\delta$ be a timelike rotational Weingarten surfaces in $\mathbb{S}^3_1$ defined by \eqref{spherical1}, \eqref{hyper1kind1} and 
\eqref{parabolic1} with the principal curvatures $\kappa$ and $\lambda$
satisfying $\kappa=\frac{a}{\lambda}$ and 
\begin{equation}
    \Psi_{-1}(\lambda,\delta)=\psi_0\pm 
    c^2\int{\frac{\lambda^2 d\lambda}{\sqrt{\lambda^2+a}(c^2-\delta(\lambda^2+a))
    \sqrt{(\lambda^2+1)c^2-\delta(\lambda^2+a)}}}
\end{equation}
for constants $\psi_0$, $a$ and $c>0$.
Then, the coordinate functions of the profile curve of ${\bf M}_\delta$ with 
constant Gaussian curvature in $\mathbb{S}^3_1$ parameterized with respect to 
$\lambda$ are given by one of the following:
\begin{itemize}
    \item [1.] for $\delta=1$,
    \subitem (a.) $y(\lambda)=\frac{c}{\sqrt{\lambda^2+a}},\;\;
        z(\lambda)=\sqrt{1-y^2(\lambda)}
        \cosh{\Psi_{-1}(\lambda,1)}\;\;
        w(\lambda)=\sqrt{1-y^2(\lambda)}
        \sinh{\Psi_{-1}(\lambda,1)}$
        where
     \subsubitem i) $\lambda<-\sqrt{c^2-a}$ or $\lambda>\sqrt{c^2-a}$ for $c^2>\mbox{max}\{1,a\}$,
     
     \subsubitem ii) $-\sqrt{\frac{c^2-a}{1-c^2}}<\lambda<-\sqrt{c^2-a}$ or 
     $\sqrt{c^2-a}<\lambda<\sqrt{\frac{c^2-a}{1-c^2}}$ for
     $a<c^2<1$,
     
     \subsubitem iii) 
     $\lambda<-\sqrt{\frac{a-c^2}{c^2-1}}$ or 
     $\lambda>\sqrt{\frac{a-c^2}{c^2-1}}$ if 
     $1<c^2<a$,
    
    \subitem (b.) 
     $y(\lambda)=\frac{c}{\sqrt{\lambda^2+a}},\;\; 
        z(\lambda)=\sqrt{y^2(\lambda)-1}
        \sinh{\Psi_{-1}(\lambda,1)},\;\mbox{and}\;
        w(\lambda)=\sqrt{y^2(\lambda,-1)}
        \cosh{\Psi_{-1}(\lambda,1)}$
        where
        \subsubitem i) $-\sqrt{c^2-a}<\lambda<\sqrt{c^2-a}$ for
        $0\leq a<c^2$,
        \subsubitem (b) $-\sqrt{c^2-a}<\lambda<-\sqrt{-a}$
        or $\sqrt{-a}<\lambda<\sqrt{c^2-a}$
        for $a<0$ and $c<1$,
        
    \item [2.] for $\delta=-1$,
        $$w(\lambda)=\frac{c}{\sqrt{\lambda^2+a}},\;\; 
        x(\lambda)=\sqrt{1+w^2(\lambda)}
        \cos{\Psi_{-1}(\lambda,-1)},\;\mbox{and}\;\;
        y(\lambda)=\sqrt{1+w^2(\lambda)}
        \sin{\Psi_{-1}(\lambda,-1)}$$
       where 
       \subitem i) $\lambda<-\sqrt{-a}$ or $\lambda>\sqrt{-a}$ for $a<0$,
       
       \subitem ii) $\lambda\in\mathbb{R}$ for $a\geq 0$,
        
     \item [3.] for $\delta=0$,
     $$z(\lambda)=\frac{c}{\sqrt{\lambda^2+a}},\;
        x(\lambda)=\Psi_{-1}(\lambda,0)z(\lambda),
        \;\;
        w(\lambda)=\frac{1-x^2(\lambda)}{2z(\lambda)}$$
        where
        \subitem i) $\lambda<-\sqrt{-a}$ or $\lambda>\sqrt{-a}$ for $a<0$,
        
        \subitem ii) $\lambda\in\mathbb{R}$ for $a\geq 0$.
     \end{itemize}
\end{Corollary}
Similarly taking account of $f(\lambda)=a\lambda^m$ for the hyperbolic rotational surface of second kind $M$ in $\mathbb{S}^3_1$ defined by \eqref{hyper2kind1}, 
the equation \eqref{hyper2kind2} becomes the following equation
\begin{equation}
    w(\lambda)=c(\lambda^{1-m}+a)^{\frac{1}{m-1}}
\end{equation}
with a positive constant $c$ and $m\neq -1$. 
Let say $\bar{S}(\lambda)=\lambda^{1-m}+a$. 
Then, we write $w(\lambda)=c\bar{S}(\lambda)^{1/{m-1}}$. 
Also, the function $\phi(\lambda)$ in \eqref{hyper2kind4} gives as
\begin{equation}
    \label{hyper2kind6}
    \phi(\lambda)=\phi_0\pm c^2
    \int{\frac{\lambda^{1-m} \bar{S}(\lambda)^{\frac{3-m}{m-1}} d\lambda}
    {(1-c^2\bar{S}(\lambda)^{2/{m-1}})
    \sqrt{1-c^2(1+\lambda^2)\bar{S}(\lambda)^{2/{(m-1)}}}}}.
\end{equation}
The remaining coordinate functions of the profile curve
are given by \eqref{hyper2kind3}.

Again, we will take $m=-1$. 
Using Theorem \ref{Weingarten1}, we find the components of profile curve of 
the hyperbolic rotational surface of second kind with constant Gaussian curvature in $\mathbb{S}^3_1$ parameterized by the principal curvature $\lambda$ seen as below. 

\begin{Corollary}
Let $M$ be a Weingarten hyperbolic rotational of second kind in a 3--dimensional de Sitter space $\mathbb{S}^3_1$ defined by \eqref{hyper2kind1} with the principal curvatures $\kappa$ and $\lambda$ satisfying $\kappa=\frac{a}{\lambda}$
with constant $a$. Then, the coordinate functions of the profile curve
of $M$ with constant Gaussian curvature in $\mathbb{S}^3_1$ parameterized with respect to $\lambda$ are given as follows:
    \begin{equation}
        w(\lambda)=\frac{c}{\sqrt{\lambda^2+a}},\;\; 
        x(\lambda)=\sqrt{1-w^2(\lambda)}\cos{\phi(\lambda)}\;\;
        y(\lambda)=\sqrt{1-w^2(\lambda)}\sin{\phi(\lambda)}
    \end{equation}
where 
\begin{equation}
    \phi(\lambda)=\phi_0\pm\frac{c^2}{\sqrt{2}}
    \int{
    \frac{\lambda^2 d\lambda}
    {(\lambda^2+a-c^2)\sqrt{\lambda^2+a}\sqrt{(1-c^2)\lambda^2+a-c^2}}}
\end{equation}
for constants $\phi_0$ and $c>0$ with 
\subitem i) 
$\displaystyle{\lambda<-\sqrt{\frac{c^2-a}{1-c^2}}}$ or 
$\displaystyle{\lambda>\sqrt{\frac{c^2-a}{1-c^2}}}$ for $a<c^2<1$ with $a<1$,

\subitem ii) $-\sqrt{\frac{a-c^2}{c^2-1}}<\lambda
<\sqrt{\frac{a-c^2}{c^2-1}}$ for $1<c^2<a$ with $a>1$,

\subitem ii) $\lambda\in\mathbb{R}$ if $c^2\leq\mbox{min}\{a,1\}$.
\end{Corollary}

\textbf{Acknowledgements.} The author would like to thank
U\u{g}ur Dursun for giving insightful comments which help to improve the quality of this manuscript.

\end{document}